\documentclass[preprint,12pt]{elsarticle}

\newtheorem{theorem}{\textbf{Theorem}}

\newtheorem{conjecture}{\textbf{Conjecture}}
\newtheorem{lemma}{\textbf{Lemma}}

\usepackage{amsfonts}
\usepackage{amsmath}

\usepackage{mathrsfs}

\def\a {\alpha}
\def\p {\phi}

\def\m {\mu}
\def\g {\gamma}
\def\N {\mathbb{N}}

\def\Q {\mathbb{Q}}

\def\e {\epsilon}

\def\g {\gamma}

\def\dfrac {\displaystyle\frac }

\usepackage{amssymb}

\journal{Bulletin of the Brazilian Mathematical Society}

\begin{document}

\begin{frontmatter}



\title{The proof of a conjecture concerning the intersection of $k$-generalized Fibonacci sequences}


\author{Diego Marques\fnref{dmm}}{\ead{diego@mat.unb.br}}
\address{Departamento de Matem\' atica, Universidade  de Bras\' ilia, Bras\' ilia, 70910-900, Brazil}

\fntext[dmm]{Supported by FAP-DF, FEMAT and CNPq-Brazil}



\begin{abstract}
For $k\geq 2$, the $k$-generalized Fibonacci sequence $(F_n^{(k)})_{n}$ is defined by the initial values $0,0,\ldots,0,1$ ($k$ terms) and such that each term afterwards is the sum of the $k$ preceding terms. In 2005, Noe and Post conjectured that the only solutions of Diophantine equation $F_m^{(k)}=F_n^{(\ell)}$, with $\ell>k>1, n>\ell+1,\ m>k+1$ are
\[
(m,n,\ell,k)=(7,6,3,2)\ \mbox{and}\ (12,11,7,3).
\]
In this paper, we confirm this conjecture.
\end{abstract}


\end{frontmatter}

\section{Introduction}\label{sec1}

Let $k\geq 2$ and denote $F^{(k)}:=(F_n^{(k)})_{n\geq -(k-2)}$, the \textit{$k$-generalized Fibonacci sequence} whose terms satisfy the recurrence relation
\begin{equation}\label{rec}
F_{n+k}^{(k)}=F_{n+k-1}^{(k)}+F_{n+k-2}^{(k)}+\cdots + F_{n}^{(k)},
\end{equation}
with initial conditions $0,0,\ldots,0,1$ ($k$ terms) and such that the first nonzero term is $F_1^{(k)}=1$.

The above sequence is one among the several generalizations of Fibonacci numbers. Such a sequence is also called $k$-\textit{step Fibonacci sequence}, the \textit{Fibonacci $k$-sequence}, or $k$-\textit{bonacci sequence}. Clearly for $k=2$, we obtain the well-known Fibonacci numbers $F_n^{(2)}=F_n$, and for $k=3$, the Tribonacci numbers $F_n^{(3)}=T_n$.

Several authors have worked on problems involving $k$-generalized Fibonacci sequences. For instance, Togb\' e and the author \cite{MT} proved that only finitely many terms of a linear recurrence sequence whose characteristic polynomial has a simple positive dominant root can be \textit{repdigits} (i.e., numbers with only one distinct digit in its decimal expansion). As an application, since the characteristic polynomial of the recurrence in \eqref{rec}, namely $x^k-x^{k-1}-\cdots - x-1$, has just one root $\alpha$ such that $|\a| > 1$ (see for instance \cite{wolf}), then there exist only finitely many terms of $F^{(k)}$ which are repdigits, for all $k\geq 2$. F. Luca \cite{Luca} and the author \cite{util} proved that $55$ and $44$ are the largest repdigits in the sequences $F^{(2)}$ and $F^{(3)}$, respectively. Moreover, the author conjectured that there are no repdigits, with at least two digits, belonging to $F^{(k)}$, for $k>3$. In a recent work, Bravo and Luca \cite{BL} confirmed this conjecture. 

Here, we are interested in the problem of determining the intersection of two $k$-generalized Fibonacci sequences. It is important to notice that Mignotte (see \cite{mig}) showed that if $(u_n)_{n\geq 0}$ and $(v_n)_{n\geq 0}$ are two linearly recurrence sequences then, under some weak technical assumptions, the equation
\begin{center}
$u_n = v_m$ 
\end{center}
has only finitely many solutions in positive integers $m,n$. Moreover, all such solutions are effectively computable (we refer the reader to \cite{A,rf16,rf17,rf15} for results on the intersection of two recurrence sequences). Thus, it is reasonable to think that the intersection $F^{(k)}\cap F^{(\ell)}$ is a finite set for all $2\leq k< \ell$. In 2005, Noe and Post \cite{noe} gave a heuristic argument to show that the expected cardinality of this intersection must be small. Furthermore, they raised the following conjecture
\begin{conjecture}[Noe-Post]\label{NP}
The Diophantine equation
\begin{equation}\label{main}
F_m^{(k)}=F_n^{(\ell)},
\end{equation}
with $\ell>k\geq 2$, $n>\ell+1$ and $m>k+1$, has only the solutions:
\begin{equation}\label{sol}
(m,n,\ell,k)=(7,6,3,2)\ \mbox{and}\ (12,11,7,3).
\end{equation}
That is,
\begin{center}
$13=F_7^{(2)}=F_6^{(3)}$\ \ and\ \ $504=F_{12}^{(3)}=F_{11}^{(7)}$
\end{center}
\end{conjecture}

Since the first nonzero terms of $F^{(k)}$ are $1,1,2,\ldots,2^{k-1}$, then the above conjecture can be rephrased as
\begin{conjecture}
Let $2\leq k<\ell$ be positive integer numbers. Then
$$
F^{(k)}\cap F^{(\ell)} = \left\{
\begin{array}{rcl}
\{0,1,2,13\},& \mbox{if} & (k,\ell)=(2,3)\\
\{0,1,2,4,504\},& \mbox{if} & (k,\ell)=(3,7)\\
\{0,1,2,8\}, & \mbox{if} & k=2\ \mbox{and}\ \ell>3\\
\{0,1,2,\ldots,2^{k-1}\}, & otherwise
\end{array}
\right.
$$
\end{conjecture}

We remark that this intersection was confirmed for $(k,\ell)=(2,3)$, by the author \cite{bohe}. Also, Noe and Post used computational methods to study this intersection (see Section \ref{program}). For our purpose we need a variant of their result which will be stated as a lemma, since we shall use it throughout our work.
\begin{lemma}\label{post} 
The only solutions $(m,n,\ell,k)$ in positive integers of Diophantine equation \eqref{main}, with $\ell>k>1, n>\ell+1,\ m>k+1, \max\{\ell,k\}<5000$ and $\max\{m,n\}<2^{5000}$, are listed in \eqref{sol}.
\end{lemma}

In this paper, we shall use transcendental tools to prove the Noe-Post conjecture. For the sake of preciseness, we stated it as a theorem.
\begin{theorem}\label{THM}
Conjecture \ref{NP} is true.
\end{theorem}

Let us give a brief overview of our strategy for proving Theorem \ref{THM}. First, we use a Dresden formula \cite[Formula (2)]{dres} to get an upper bound for a linear form in three logarithms related to equation \eqref{main}. After, we use a lower bound due to Matveev to obtain an upper bound for $m$ and $n$ in terms of $\ell$. Very recently, Bravo and Luca solved the equation $F_n^{(k)}=2^m$ and for that they used a nice argument combining some estimates together with the Mean Value Theorem (this can be seen in pages 72 and 73 of \cite{BL2}). In our case, we must use two times this Bravo and Luca approach together with a reduction argument due to Dujella and Peth$\ddot{\mbox{o}}$ to prove our main theorem. In the final section, we present a program for checking the ``small"\ cases. The computations in the paper were performed using \textit{Mathematica}$^{\tiny {\textregistered}}$.

We remark some differences between our work and the one by Bravo and Luca. In their paper, the equation $F_n^{(k)}=2^m$ was studied. By applying a key method, they get directly an upper bound for $|2^m-2^{n-2}|$. In our case, the equation $F_m^{(k)}=F_n^{(\ell)}$ needs a little more work, because it is necessary to apply two times their method to get an upper bound for $|2^{n-2}-2^{m-2}|$. Moreover, they used a reduction argument due to Dujella and Peth$\ddot{\mbox{o}}$ to solve all small cases. In our work, we use a Noe and Post program to deal with the ``very"\ small cases. Our presentation is therefore organized in a similar way that the one in the papers \cite{BL2,BL}, since we think that those presentations are intuitively clear.

\section{Upper bounds for $m$ and $n$ in terms of $\ell$}\label{sec2}

In this section, we shall prove the following result
\begin{lemma}\label{t1}
If $(m,n,\ell,k)$ is a solution in positive integers of Diophantine equation \eqref{main}, with $\ell>k\geq 2$, $n>\ell+1$ and $m>k+1$. Then
\begin{equation*}
n<m<4.4\cdot 10^{14}\ell^8\log^3 \ell.
\end{equation*}
\end{lemma}

Before proceeding further, we shall recall some facts and properties of these sequences which will be used
after.

We know that the characteristic polynomial of $(F_n^{(k)})_n$ is 
\[
\psi_k(x):=x^k-x^{k-1}-\cdots -x-1
\]
and it is irreducible over $\Q[x]$ with just one zero outside the unit circle. That single zero is located between $2(1-2^{-k})$ and $2$ (as can be seen in \cite{wolf}). Also, in a recent paper, G. Dresden \cite[Theorem 1]{dres} gave a simplified ``Binet-like"\ formula for $F_n^{(k)}$:
\begin{equation}\label{binet}
F_n^{(k)}=\displaystyle\sum_{i=1}^k\dfrac{\alpha_i-1}{2+(k+1)(\a_i-2)}\a_i^{n-1},
\end{equation}
for $\alpha=\alpha_1,\ldots,\a_k$ being the roots of $\psi_k(x)$. There are many other ways of
representing these $k$-generalized Fibonacci numbers, as can be seen in \cite{2,3,4,5}. Also, it was proved in \cite[Lemma 1]{BL} that
\begin{equation}\label{Lu}
\alpha^{n-2}\leq F_n^{(k)}\leq \alpha^{n-1},\ \mbox{for\ all}\ n\geq 1,
\end{equation}
where $\a$ is the dominant root of $\psi_k(x)$. Also, the contribution of the roots inside the unit circle in formula \eqref{binet} is almost trivial. More precisely, it was proved in \cite{dres} that
\begin{equation}\label{small}
|F_n^{(k)}-g(\alpha,k)\alpha^{n-1}|<\dfrac{1}{2},
\end{equation}
where we adopt throughout the notation $g(x,y):=(x-1)/(2+(y+1)(x-2))$.

Another tool to prove Lemma \ref{t1}, we still use a lower bound for a linear form logarithms {\it \`a la Baker} and such a bound was given by the following result of Matveev (see \cite{matveev} or Theorem 9.4 in \cite{bugeaud}).

\begin{lemma}\label{lemma1}
Let $\gamma_1,\ldots,\gamma_t$ be real algebraic numbers and let $b_1,\ldots,b_t$ be nonzero rational integer numbers. Let $D$ be the degree of the number field $\Q(\gamma_1,\ldots,\gamma_t)$ over $\Q$ and let $A_j$ be a positive real number satisfying
\begin{center}
$A_j\geq \max\{Dh(\g_j),|\log \g_j|,0.16\}$ for $j=1,\ldots,t$.
\end{center}
Assume that
$$
B\geq \max\{|b_1|,\ldots,|b_t|\}.
$$
If $\g_1^{b_1}\cdots \g_t^{b_t}\neq 1$, then
$$
|\g_1^{b_1}\cdots \g_t^{b_t}-1|\geq \exp(-1.4\cdot 30^{t+3}\cdot t^{4.5}\cdot D^2(1+\log D)(1+\log B)A_1\cdots A_t).
$$
\end{lemma}
As usual, in the above statement, the \textit{logarithmic height} of an $s$-degree algebraic number $\g$ is defined as
$$
h(\g)=\dfrac{1}{s}(\log |a|+\displaystyle\sum_{j=1}^s\log \max\{1,|\g^{(j)}|\}),
$$
where $a$ is the leading coefficient of the minimal polynomial of $\gamma$ (over $\mathbb{Z}$) and $(\g^{(j)})_{1\leq j\leq s}$ are the conjugates of $\g$ (over $\Q$).

\subsection{The proof of Lemma \ref{t1}}

First, the inequality $n<m$ follows from the facts that the sequences $(F_n^{(\ell)})_n$ and $(F_n^{(\ell)})_{\ell}$ are nondecreasing together with \eqref{main}, $n>\ell+1$ and $m>k+1$. By the way, to find an upper bound for $m$ in terms of $n$, we combine \eqref{main} and \eqref{Lu} to obtain
\begin{equation}\label{ii}
2^{n-1}>\p^{n-1}\geq F_n^{(\ell)}=F_m^{(k)}\geq \a^{m-2}>(\sqrt{2})^{m-2}\ \mbox{and\ so}\ 2n>m,
\end{equation}
where in the last inequality we used that $\a>3/2>\sqrt{2}$.

Now, we use \eqref{small} to get
\begin{center}
$|F_m^{(k)}-g(\alpha,k)\alpha^{m-1}|<\dfrac{1}{2}$ and $|F_n^{(\ell)}-g(\phi,\ell)\p^{n-1}|<\dfrac{1}{2}$,
\end{center}
where $\a$ and $\p$ are the dominant roots of the recurrences $(F_m^{(k)})_m$ and $(F_n^{(\ell)})_n$, respectively. Combining these inequalities, we obtain
\begin{equation}\label{11}
|g(\p,\ell)\p^{n-1}-g(\alpha,k)\alpha^{m-1}|<1
\end{equation}
and so
\begin{equation}\label{12}
\left|\frac{g(\p,\ell)\p^{n-1}}{g(\alpha,k)\alpha^{m-1}}-1\right |<\frac{1}{g(\alpha,k)\alpha^{m-1}}<\frac{4}{\alpha^{m-1}},
\end{equation}
where we used that $g(\alpha,k)>1/4$, since $\alpha> 3/2$ (for $k\geq 2$) and $2+(k+1)(\a-2)<2$. Thus \eqref{12} becomes
\begin{equation}\label{upper}
|e^{\Lambda}-1|<\frac{4}{\a^{m-1}},
\end{equation}
where $\Lambda:=(n-1)\log \p+\log(g(\p,\ell)/g(\a,k))-(m-1)\log \a$. 

Now, we shall apply Lemma \ref{lemma1}. To this end, take $t:=3$,
\[
\g_1:=\p,\ \g_2:=\dfrac{g(\p,\ell)}{g(\a,k)},\ \g_3:=\a
\]
and 
\[
b_1:=n-1,\ b_2:= 1,\ b_3:= m-1.
\]
For this choice, we have $D=[\Q(\a,\p):\Q]\leq k\ell<\ell^2$. Also $h(\g_1)=(\log \p)/\ell<(\log 2)/\ell<0.7/\ell$ and similarly $h(\g_3)<0.7/k$. In \cite[p. 73]{BL2}, an estimate for $h(g(\a,k))$ was given. More precisely, it was proved that
\[
h(g(\a,k)) < \log (k+1)+\log 4.
\]
Analogously, 
\[
h(g(\p,\ell))<\log (\ell+1)+\log 4.
\]
Thus
\[
h(\g_2)\leq h(g(\p,\ell))+h(g(\a,k))\leq \log (\ell+1)+\log (k+1)+2\log 4,
\]
where we used the well-known facts that $h(xy)\leq h(x)+h(y)$ and $h(x)=h(x^{-1})$. Also, in \cite{BL2} was proved that $|g(\a_i,k)|<2$, for all $i=1,\ldots, k$.

Since $\ell>k$ and $m>n$, we can take $A_1=A_3:= 0.7\ell,\ A_2:=2\ell^2\log (4\ell+4)$ and $B:=m-1$. 

Before applying Lemma \ref{lemma1}, it remains us to prove that $e^{\Lambda}\neq 1$. Suppose, towards a contradiction, the contrary, i.e., $g(\a,k)\a^{m-1}=g(\p,\ell)\p^{n-1}\in \Q(\p)$. So, we can conjugate this relation in $\Q(\p)$ to get
\begin{center}
$g(\a_{s_i},k)\a_{s_i}^{m-1}=g(\p_i,\ell)\p_i^{n-1}$, for $i=1,\ldots,\ell$,
\end{center}
where $\a_{s_i}$ are the $\ell$ conjugates of $\a$ over $\Q(\p)$. Since $g(\a,k)\a^{m-1}$ has at most $k$ conjugates (over $\Q$), then each number in the list $\{g(\a_{s_i},k)\a_{s_i}^{m-1}:1\leq i\leq \ell\}$ is repeated at least $\ell/k>1$ times. In particular, there exists $t\in \{2,\ldots,\ell\}$, such that $g(\a_{s_1},k)\a_{s_1}^{m-1}=g(\a_{s_t},k)\a_{s_t}^{m-1}$. Thus, $g(\p,k)\p^{n-1}=g(\p_t,\ell)\p_t^{n-1}$ and then
\[
\left(\dfrac{7}{4}\right)^{n-1}<\p^{n-1}=\left|\dfrac{g(\p_t,\ell)}{g(\p,\ell)}\right||\p_t|^{n-1}<8,
\]
where we used that $\p>2(1-2^{-\ell})\geq 7/4,\ |g(\p_t,\ell)|<2<8|g(\p,\ell)|$ and $|\p_t|<1$ for $t>1$. However, the inequality $(7/4)^{n-1}<8$ holds only for $n=1,2,3,4$, but this gives an absurdity, since $n>\ell+1\geq 3+1=4$. Therefore $e^{\Lambda}\neq 1$.

Now, the conditions to apply Lemma \ref{lemma1} are fulfilled and hence
\[
|e^{\Lambda}-1|>\exp(-1.5\cdot 10^{11}\ell^8(1+2\log \ell)\log (4\ell+4)(1+\log (m-1)))
\]
Since, $1+2\log \ell\leq 3\log \ell$ and $4\ell+4<\ell^{2.6}$ (for $\ell\geq 3$), we have that 
\begin{equation}\label{lower}
|e^{\Lambda}-1|>\exp(-2.4\cdot 10^{12}\ell^8\log^2 \ell \log (m-1))
\end{equation}
By combining \eqref{upper} and \eqref{lower}, we get
\[
\dfrac{m-1}{\log (m-1)}<6.1\cdot 10^{12}\ell^8\log^2 \ell,
\]
where we used that $\log \alpha>0.4$. Since the function $x/\log x$ is increasing for $x>e$, it is a simple matter to prove that
\begin{equation}\label{key}
\dfrac{x}{\log x}<A\ \ \mbox{implies\ that}\ \ x<2A\log A.
\end{equation}
A proof for that can be found in \cite[p. 74]{BL2}.

Thus, by using \eqref{key} for $x:=m-1$ and $A:=6.1\cdot 10^{12}\ell^8\log^2 \ell$, we have that
\[
m-1<2(6.1\cdot 10^{12}\ell^8\log^2 \ell)\log (6.1\cdot 10^{12}\ell^8\log^2 \ell).
\]
Now, the inequality $30+2\log \log \ell<28\log \ell$, for $\ell\geq 3$, yields
\[
\log (6.1\cdot 10^{12}\ell^8\log^2 \ell)<30+8\log \ell+2\log \log \ell<36\log \ell.
\]
Therefore
\begin{equation}\label{m<l}
m<4.4\cdot 10^{14}\ell^8\log^3 \ell
\end{equation}
The proof is then complete.
\qed

\section{Upper bound for $\ell$ in terms of $k$}\label{sec3}

\begin{lemma}\label{t3}
If $(m,n,\ell,k)$ is a solution in positive integers of equation \eqref{main}, with $\ell>k>1, n>\ell+1$ and $m>k+1$, then
\begin{equation}\label{l<k}
\ell<1.8\cdot 10^{16}k^3\log^3 k.
\end{equation}
\end{lemma}

\noindent
{\bf Proof.} If $\ell\leq 239$, then the inequalities \eqref{m<l} yields $m<8\cdot 10^{35}$. Since $k<\ell$ and $n<m$, then Lemma \ref{post} ensures that the only solutions of equation \eqref{main} with the conditions in the statement of Theorem \ref{THM} are $(m,n,\ell,k)=(7,6,3,2)$ and $(12,11,7,3)$. 

Thus, we may assume that $\ell>239$. Therefore
\begin{equation}\label{beg}
n<4.4\cdot 10^{14}\ell^8\log^3 \ell<2^{\ell/2}
\end{equation}
where we used \eqref{m<l} and the fact that $n<m$. By using a key argument due to Bravo and Luca \cite[p. 72-73]{BL2}, we get
\begin{equation}\label{Phi}
|2^{n-2}-g(\a,k)\a^{m-1}|<\dfrac{5\cdot 2^{n-2}}{2^{\ell/2}}
\end{equation}
or equivalently,
\begin{equation}\label{U}
|1-g(\a,k)\a^{m-1}2^{-(n-2)}|<\dfrac{5}{2^{\ell/2}}.
\end{equation}

For applying Lemma \ref{lemma1}, it remains us to prove that the left-hand side of \eqref{U} is nonzero, or equivalently, $2^{n-2}\neq g(\a,k)\a^{m-1}$. To obtain a contradiction, we suppose the contrary, i.e., $2^{n-2}=g(\a,k)\a^{m-1}$. By conjugating the previous relation in the splitting field of $\psi_k(x)$, we obtain $2^{n-2}=g(\a_i,k)\a_i^{m-1}$, for $i=1,\ldots,k$. However, when $i>1$, $|\a_i|<1$ and $|g(\a_i,k)|<2$. But this leads to the following absurdity
\[
2^{n-2}=|g(\a_i,k)||\a_i|^{m-1}<2,
\]
since $n> 4$. Therefore $g(\a,k)\a^{m-1}2^{-(n-2)}\neq 1$ and then we are in position to apply Lemma \ref{lemma1}. For that, take $t:=3$,
\[
\g_1:=g(\a,k),\ \g_2:= \a,\ \g_3:=2
\]
and
\[
b_1:= 1,\ b_2:=m-1,\ b_3:=-(n-2).
\]

By some calculations made in Section \ref{sec2}, we see that $A_1:=k\log (4k+4),\ A_2=A_3:=0.7$ are suitable choices. Moreover $D= k$ and $B=m-1$. Thus
\begin{equation}\label{L}
|1-g(\a,k)\a^{m-1}2^{-(n-2)}|>\exp(-C_1k^3(1+\log k)(1+\log (m-1))\log (4k+4)),
\end{equation}
where we can take $C_1=0.75\cdot 10^{11}$. Combining \eqref{U} and \eqref{L} together with a straightforward calculation, we get
\begin{equation}\label{back}
\ell<4.7\cdot 10^{12}k^3\log^2 k\log m
\end{equation}

On the other hand, $m<4.4\cdot 10^{14}\ell^8 \log^3 \ell$ (by \eqref{m<l}) and so
\begin{equation}\label{ss}
\log m<\log (4.4\cdot 10^{14}\ell^8 \log^3 \ell)<45\log \ell.
\end{equation}
Turning back to inequality \eqref{back}, we obtain
\[
\dfrac{\ell}{\log \ell}<2.2\cdot 10^{14}k^3\log^2 k
\]
which implies (by \eqref{key}) that 
\[
\ell<2(2.2\cdot 10^{14}k^3\log^2 k)\log (2.2\cdot 10^{14}k^3\log^2 k).
\]
Since $\log (2.2\cdot 10^{14}k^3\log^2 k)<39\log k$, we finally get the desired inequality
\[
\ell<1.8\cdot 10^{16}k^3\log^3 k.
\]
\qed

\section{The proof of Theorem \ref{THM}}\label{sec4}

In order to finish the proof of Theorem \ref{THM}, our last ingredient is a variant of the famous Baker-Davenport lemma, which is due to Dujella and Peth\H{o} \cite{dujella}. For a real number $x$, we use  $\parallel x \parallel= \min\{|x-n|:n\in \N\}$ for the distance from $x$ to the nearest integer.
\begin{lemma}\label{lemma2}
Suppose that $M$ is a positive integer. Let $p/q$ be a convergent of the continued fraction expansion of the irrational number $\gamma$ such that $q > 6M$ and let $A, B$ be some real numbers with $A>0$ and $B>1$. Let $\epsilon=\parallel \mu q \parallel-M\parallel \gamma q \parallel$, where $\mu$ is a real number. If $\epsilon>0$, then there is no solution to the inequality
$$
0<m\gamma -n+\mu < A\cdot B^{-k}
$$
in positive integers $m,n$ and $k$ with
\begin{center}
$m\leq M$ and $k\geq \dfrac{\log(Aq/\epsilon)}{\log B}$.
\end{center}
\end{lemma}
See Lemma 5, a.) in \cite{dujella}.

The proof of our main result splits in two cases:

\subsection{The case $k>1655$}

First, let us prove that there is no solution when $k>1655$. Towards a contradiction, suppose that $(m,n,\ell,k)$ is such a solution. Then the inequality $\ell<1.8\cdot 10^{16}k^3\log^3 k$ together with \eqref{m<l} yield
\begin{eqnarray*}
m & < & 4.4\cdot 10^{14} (1.8\cdot 10^{16}k^3\log^3 k)^8\log^3 (1.8\cdot 10^{16}k^3\log^3 k)\\
 & < & 3\cdot 10^{148}k^{24}\log^{27} k<2^{k/2},
\end{eqnarray*}
where the last inequality holds only because $k>1655$. Now, we use again the key argument of Bravo and Luca to conclude that
\begin{equation}\label{Alf}
|2^{m-2}-g(\p,\ell)\p^{n-1}|<\dfrac{5\cdot 2^{m-2}}{2^{k/2}}.
\end{equation}
Combining \eqref{Phi}, \eqref{Alf} and \eqref{11}, we get
\begin{eqnarray*}
|2^{n-2}-2^{m-2}| & \leq & |2^{n-2}-g(\a,k)\a^{n-1}|+|g(\a,k)\a^{n-1}-g(\p,\ell)\p^{n-1}|\\
& & + |2^{m-2}-g(\p,\ell)\p^{n-1}|\\
& < & \dfrac{5\cdot 2^{n-2}}{2^{\ell /2}} + 1 + \dfrac{5\cdot 2^{m-2}}{2^{k/2}} < \dfrac{11\cdot 2^{m-2}}{2^{k/2}},
\end{eqnarray*}
since $n<m,\ k<\ell$ and $m>k+1$. Therefore
\begin{equation}\label{m-n}
|2^{n-m}-1|<\dfrac{11}{2^{k/2}}.
\end{equation}

Since $n\leq m-1$, then
\[
\dfrac{1}{2}\leq 1-2^{n-m}=|2^{n-m}-1|<\dfrac{11}{2^{k/2}}.
\]
Thus $2^{k/2}<22$ leading to an absurdity, since $k>1655$. 

\subsection{The case $2\leq k\leq 1655$}

If $k\leq 1655$, then $\ell<4\cdot 10^{28}$ (by \eqref{l<k}). Thus, by \eqref{m<l}, one has that $n<m<2\cdot 10^{248}$. In order to use the Lemma \ref{lemma2}, we rewrite (\ref{U}) as
\[
|e^{\Theta}-1|<\dfrac{5}{2^{\ell/2}},
\]
where $\Theta:=(m-1)\log \a-(n-2)\log 2+\log g(\a,k)$. Recall that we proved that $e^{\Theta}\neq 1$ (the paragraph below (\ref{U})) and so $\Theta\neq 0$. 

If $\Theta>0$, then $\Theta<e^{\Theta}-1<5/2^{\ell/2}$. In the case of $\Theta<0$, we use $1-e^{-|\Theta|}=|e^{\Theta}-1|<5/2^{\ell/2}$ to get $e^{|\Theta|}<1/(1-5\cdot 2^{-\ell/2})$. Thus
\[
|\Theta|<e^{|\Theta|}-1<\dfrac{5\cdot 2^{-\ell/2}}{1-5\cdot 2^{-\ell/2}}<5\cdot 2^{-\ell/2+4},
\]
where we used that $1/(1-5\cdot 2^{-\ell/2})<16$, for $\ell\geq 3$. Summarizing, the further arguments work for $\Theta>0$ and $\Theta<0$ in a very similar way.

Thus, to avoid unnecessary repetitions we shall consider only the case $\Theta>0$. For that, we have
\[
0<(m-1)\log \a-(n-2)\log 2+\log g(\a,k)< 5\cdot (\sqrt{2})^{-\ell}
\]
and then
\begin{equation}\label{A}
0<(m-1)\g_k-(n-2)+\m_k<7.3\cdot (\sqrt{2})^{-\ell},
\end{equation}
with $\g_k:=\log \a^{(k)}/\log 2$ and $\m_k:=\log g(\a^{(k)},k)/\log 2$. Here, we added the superscript to $\a$ for emphasizing its dependence on $k$.

We claim that $\g_k$ is irrational, for any integer $k\geq 2$. In fact, if $\g_k=p/q$, for some positive integers $p$ and $q$, we have that $2^p=(\a^{(k)})^q$ and as before we can conjugate this relation by some automorphism of the Galois group of the splitting field of $\psi_k(x)$ over $\Q$ to get $2^p<|(\a_i^{(k)})^q|<1$, for $i>1$, which is an absurdity, since $p\geq 1$. Let $q_{n,k}$ be the denominator of the $n$-th convergent of the continued fraction of $\g_k$. Taking $M_k:=3\cdot 10^{148}k^{24}\log^{27} k\leq M_{1655}<2\cdot 10^{248}$, we use \textit{Mathematica} to get
\[
\displaystyle\min_{2\leq k\leq 1655}q_{650,k}>6\cdot 10^{308}>6M_{1655}.
\]
Also
\[
\displaystyle\max_{2\leq k\leq 1655}q_{650,k}<2\cdot 10^{1125}.
\]

Define $\e_k:=\parallel \m_kq_{650,k}\parallel-M_k\parallel \g_kq_{650,k} \parallel$, for $2\leq k\leq 1655$, we get
\[
\displaystyle\min_{2\leq k\leq 1655}\e_{k}>0.000015.
\]

Note that the conditions to apply Lemma \ref{lemma2} are fulfilled for $A=7.3$ and $B=\sqrt{2}$, and hence there is no solution to inequality \eqref{A} (and then no solution to the Diophantine equation \eqref{main}) for $m$ and $\ell$ satisfying
\begin{center}
$m<M_k< 2\cdot 10^{248}$ and $\ell \geq \dfrac{\log (Aq_{650,k}/\e_k)}{\log B}$.
\end{center}
Since $m<M_k$ (for $2\leq k\leq 1655$), then
\[
\ell < \dfrac{\log (Aq_{650,k}/\e_k)}{\log B}\leq \dfrac{\log (7.3\cdot 2\cdot 10^{1125}/1.5\cdot 10^{-5})}{\log \sqrt{2}}<3757.0616\ldots.
\]

Therefore $2\leq k\leq 1655$ and $\ell\leq 3757$. Now, by applying Lemma \ref{t1}, we obtain $n<m<9.8\cdot 10^{45}$. However this is a case already treated in Lemma \ref{post}. Thus, the only solutions of equation \eqref{main} with $\ell>k>1, n>\ell+1$ and $m>k+1$ are those listed in \eqref{sol}. Thus, the proof of Theorem \ref{THM} is complete.\qed

\section{The program}\label{program}

In this section, for the sake of completeness, we present the \textit{Mathematica} program (which was kindly sent to us by Noe \cite{Tony}) used to confirm Lemma \ref{post}:

\begin{verbatim}
nn = 5000;
f = 2^Range[nn] - 1;
f[[1]] = Infinity;
cnt = 0;
seq = Table[Join[2^Range[i - 1], {2^i - 1}], {i, nn}];
done = False;
While[! done, fMin = Min[f];
 pMin = Flatten[Position[f, fMin]];
 If[Length[pMin] > 1, Print[{fMin, pMin}]];
 Do[k = pMin[[i]];
  s = Plus @@ seq[[k]];
  seq[[k]] = RotateLeft[seq[[k]]];
  seq[[k, k]] = s;
  f[[k]] = s, {i, Length[pMin]}];
 cnt++;
 done = (fMin > 2^nn)]; cnt
\end{verbatim}

The calculations in this paper took roughly $132$ hours on 2.5 GHz Intel Core i5 4GB Mac OSX.

\subsection*{Acknowledgements}
We would like to express our deepest gratitude to the anonymous referee for carefully examining this paper and providing it a number of important comments. In particular, he/she pointed out the useful reference \cite{BL2}. This research was partly supported by FAP-DF, FEMAT and CNPq.


\begin{thebibliography}{HD}




\normalsize
\baselineskip=17pt


\bibitem{A} M. A. Alekseyev, \textit{On the intersections of Fibonacci, Pell, and Lucas numbers}, Integers 11A (2011).

\bibitem{BL2} J. Bravo, F. Luca, \textit{Powers of two in generalized Fibonacci sequences}, Rev. Colombiana Mat. \textbf{46} (2012), 67--79.

\bibitem{BL} J. Bravo, F. Luca, \textit{$k$-generalized Fibonacci numbers with only one distinct digit}, Preprint.

\bibitem{bugeaud} Y. Bugeaud, M. Mignotte, S. Siksek, \textit{Classical and modular approaches
to exponential Diophantine equations. I. Fibonacci and Lucas perfect powers}, Ann. of Math. \textbf{163} (2006), no. 3, 969--1018.


\bibitem{dujella} A. Dujella and A. Peth$\ddot{\mbox{o}}$, \textit{A generalization of a theorem of Baker and Davenport}, Quart. J. Math. Oxford Ser. (2) \textbf{49} (1998), 291--306. 

\bibitem{dres} G. P. Dresden, \textit{A simplified Binet formula for $k$-generalized Fibonacci numbers}, Preprint, 	arXiv:0905.0304v2 (2011). Accessed 31 December 2011.



\bibitem{2} D. E. Ferguson, \textit{An expression for generalized Fibonacci numbers}, Fibonacci Quart.
\textbf{4} (1966), 270–273.

\bibitem{3} I. Flores, \textit{Direct calculation of $k$-generalized Fibonacci numbers}, Fibonacci Quart. \textbf{5}
(1967), 259--266. 

\bibitem{4} H. Gabai, \textit{Generalized Fibonacci $k$-sequences}, Fibonacci Quart. \textbf{8} (1970), 31–38. 


\bibitem{5} D. Kalman, \textit{Generalized Fibonacci numbers by matrix methods}, Fibonacci Quart. \textbf{20}
(1982), 73–76.












\bibitem{rf7} P. Y. Lin, \textit{De Moivre-type identities for the Tribonacci numbers}. Fibonacci Quart. \textbf{26} (1988), 131-134. 

\bibitem{Luca} F. Luca, \textit{Fibonacci and Lucas numbers with only one distinct digit}, Port. Math. {\bf 57} (2) (2000), 243--254.  




\bibitem{MT} D. Marques and A. Togb\' e, \textit{On terms of a linear recurrence sequence with only one distinct block of digits}, Colloq. Math. \textbf{124}, p. 145-155, 2011.

\bibitem{bohe} D. Marques, \textit{On the intersection of two distinct $k$-generalized Fibonacci sequences}, To appear in Math. Bohem.

\bibitem{util} D. Marques, \textit{On $k$-generalized Fibonacci numbers with only one distinct digit}, To appear in Util. Math.

\bibitem{matveev} E. M. Matveev, \textit{An explicit lower bound for a homogeneous rational linear form in the logarithms of algebraic numbers}, Izv. Math. \textbf{64} (2000), no. 6, 1217--1269.

\bibitem{mig} M. Mignotte, \textit{Intersection des images de certaines suites r\' ecurrentes lin\' eaires.} Theor. Comput. Sci., \textbf{7} (1), (1978), 117-121.

\bibitem{noe} T. D. Noe and J. V. Post, \textit{Primes in Fibonacci $n$-step and Lucas $n$-step sequences}, J. Integer Seq., \textbf{8} (2005), Article 05.4.4.

\bibitem{Tony} T. D. Noe, personal communication, 27 January 2012.


\bibitem{rf16} H. P. Schlickewei, W. M. Schmidt, \textit{Linear equations in members of recurrence sequences.} 
Ann. Scuola Norm. Sup. Pisa Cl. Sci. (4) \textbf{20} (1993)2, 219-246. 

\bibitem{rf17} H. P. Schlickewei, W. M. Schmidt, \textit{The intersection of recurrence sequences.} Acta Arith. \textbf{72} (1995) 1-44. 



\bibitem{rf6} W. R. Spickerman, \textit{Binet's formula for the Tribonacci sequence}, Fibonacci Quart. \textbf{20} (1982), 118--120. 

\bibitem{rf15} S. K. Stein, \textit{The intersection of Fibonacci sequences}. Michigan Math. J. \textbf{9} (1962) 399-402. 


\bibitem{wolf} A. Wolfram, \textit{Solving generalized Fibonacci recurrences}, Fibonacci Quart. {\bf 36} (1998),
129--145. 


\end{thebibliography}
\end{document}